\documentclass[11pt]{article} 

\usepackage{amsfonts}
\usepackage{latexsym}
\catcode`\@=11\relax
\newwrite\@unused
\def\typeout#1{{\let\protect\string\immediate\write\@unused{#1}}}

\def\psglobal#1{
\immediate\special{ps: plotfile #1 }}
\def\psfiginit{\typeout{psfiginit}
\immediate\psglobal{figtex.pro}%
\special{ps:: /TeXMagnification {\the\mag} def}
}

\def\@nnil{\@nil}
\def\@empty{}
\def\@psdonoop#1\@@#2#3{}
\def\@psdo#1:=#2\do#3{\edef\@psdotmp{#2}\ifx\@psdotmp\@empty \else
    \expandafter\@psdoloop#2,\@nil,\@nil\@@#1{#3}\fi}
\def\@psdoloop#1,#2,#3\@@#4#5{\def#4{#1}\ifx #4\@nnil \else
       #5\def#4{#2}\ifx #4\@nnil \else#5\@ipsdoloop #3\@@#4{#5}\fi\fi}
\def\@ipsdoloop#1,#2\@@#3#4{\def#3{#1}\ifx #3\@nnil
       \let\@nextwhile=\@psdonoop \else
      #4\relax\let\@nextwhile=\@ipsdoloop\fi\@nextwhile#2\@@#3{#4}}
\def\@tpsdo#1:=#2\do#3{\xdef\@psdotmp{#2}\ifx\@psdotmp\@empty \else
    \@tpsdoloop#2\@nil\@nil\@@#1{#3}\fi}
\def\@tpsdoloop#1#2\@@#3#4{\def#3{#1}\ifx #3\@nnil
       \let\@nextwhile=\@psdonoop \else
      #4\relax\let\@nextwhile=\@tpsdoloop\fi\@nextwhile#2\@@#3{#4}}
\def\psdraft{
	\def\@psdraft{0}
	\def\@psdraftspecial{100}
}
\def\psdraftspecial{
	\def\@psdraft{0}
	\def\@psdraftspecial{0}
}
\def\psfull{
	\def\@psdraft{100}
}
\psfull

\newif\if@prologfile
\newif\if@postlogfile
\newif\if@bbllx
\newif\if@bblly
\newif\if@bburx
\newif\if@bbury
\newif\if@height
\newif\if@width
\newif\if@rheight
\newif\if@rwidth
\newif\if@clip
\newif\if@right
\newif\if@left
\newif\if@toplines
\newif\if@box
\newif\if@caption
\newif\if@surround
\newif\if@captionwidth
\newif\if@captionwrite
\newif\if@captionopen
\def\@p@@sclip#1{\@cliptrue}
\def\@p@@sfile#1{
		\def\@p@sfile{#1}
}
\def\@p@@sfigure#1{
		\def\@p@sfile{#1}
}
\def\@p@sfake{\hbox to 0pt{\hss Whatever\hss}}
\def\@p@@sbbllx#1{
		\@bbllxtrue
		\@d@mscratch=#1
		\edef\@p@sbbllx{\number\@d@mscratch}
}
\def\@p@@sbblly#1{
		\@bbllytrue
		\@d@mscratch=#1
		\edef\@p@sbblly{\number\@d@mscratch}
}
\def\@p@@sbburx#1{
		\@bburxtrue
		\@d@mscratch=#1
		\edef\@p@sbburx{\number\@d@mscratch}
}
\def\@p@@sbbury#1{
		\@bburytrue
		\@d@mscratch=#1
		\edef\@p@sbbury{\number\@d@mscratch}
}
\def\@p@@sheight#1{
		\@heighttrue
		\@d@mscratch=#1
   		\edef\@p@sheight{\number\@d@mscratch}
}
\def\@p@@swidth#1{
		\@widthtrue
		\@d@mscratch=#1
		\edef\@p@swidth{\number\@d@mscratch}
}
\def\@p@@srheight#1{
		\@rheighttrue
		\@d@mscratch=#1
		\edef\@p@srheight{\number\@d@mscratch}
}
\def\@p@@srwidth#1{
		\@rwidthtrue
		\@d@mscratch=#1
		\edef\@p@srwidth{\number\@d@mscratch}
}
\def\@p@@sright#1{\@righttrue \@surroundtrue}
\def\@p@@sleft#1{\@lefttrue \@surroundtrue}
\def\@p@@sextraheight#1{\@d@mextraheight=#1}
\def\@p@@sbox#1{\@boxtrue}
\def\@p@@scaption#1{\@captiontrue}
\def\@p@@stoplines#1{
		\@toplinestrue
		\@c@ttoplines=#1
}
\def\@p@@scaptionwidth#1{
		\@captionwidthtrue
	  	\@d@mcaptionwidth=#1
}
\def\@p@@scaptionwrite#1{
		\global\@captionwritetrue
		\global\@w@rname=\expandafter{\jobname_captions.tex}
		\typeout{Captions are written to \the\@w@rname}
}

\def\@p@@sprolog#1{\@prologfiletrue\def\@prologfileval{#1}}
\def\@p@@spostlog#1{\@postlogfiletrue\def\@postlogfileval{#1}}
\def\@cs@name#1{\csname #1\endcsname}
\def\@setparms#1=#2,{\@cs@name{@p@@s#1}{#2}}
\def\ps@init@parms{
		\@bbllxfalse \@bbllyfalse
		\@bburxfalse \@bburyfalse
		\@heightfalse \@widthfalse
		\@rheightfalse \@rwidthfalse
		\def\@p@sbbllx{}\def\@p@sbblly{}
		\def\@p@sbburx{}\def\@p@sbbury{}
		\def\@p@sheight{}\def\@p@swidth{}
		\def\@p@srheight{}\def\@p@srwidth{}
		\def\@p@sfile{}
		\def\@p@scost{10}
		\def\@sc{}
		\@prologfilefalse
		\@postlogfilefalse
		\@clipfalse
		\@rightfalse \@leftfalse
		\@boxfalse \@captionfalse
		\@toplinesfalse \@surroundfalse
		\@d@mextraheight=0pt
 		\@c@ttoplines=0
		\@pshape={} \def\@p@srheight@total{}
		\@captionwidthfalse \@d@mcaptionwidth=0pt
}
\def\parse@ps@parms#1{
	 	\@psdo\@psfiga:=#1\do
		   {\expandafter\@setparms\@psfiga,}}
\newif\ifno@bb
\newif\ifnot@eof
\newread\ps@stream
\newtoks\@linetok
\def\bb@missing{
	\typeout{psfig: searching \@p@sfile \space  for bounding box}
	\openin\ps@stream=\@p@sfile
	\no@bbtrue
	\not@eoftrue
	\catcode`\%=12
	\loop
		\read\ps@stream to \line@in
		\global\@linetok=\expandafter{\line@in}
		\ifeof\ps@stream \not@eoffalse \fi
		\@bbtest{\@linetok}
		\if@bbmatch\not@eoffalse\expandafter\bb@cull\the\@linetok\fi
	\ifnot@eof \repeat
	\catcode`\%=14
}	
\catcode`\%=12
\newif\if@bbmatch
\def\@bbtest#1{\expandafter\@a@\the#1
\long\def\@a@#1
     \ifx\@bbtest#2\@bbmatchfalse\else\@bbmatchtrue\fi}
\long\def\bb@cull#1 #2 #3 #4 #5 {
	\@d@mscratch=#2 bp\edef\@p@sbbllx{\number\@d@mscratch}
	\@d@mscratch=#3 bp\edef\@p@sbblly{\number\@d@mscratch}
	\@d@mscratch=#4 bp\edef\@p@sbburx{\number\@d@mscratch}
	\@d@mscratch=#5 bp\edef\@p@sbbury{\number\@d@mscratch}
	\no@bbfalse
}
\catcode`\%=14
\def\compute@bb{
		\no@bbfalse
		\if@bbllx \else \no@bbtrue \fi
		\if@bblly \else \no@bbtrue \fi
		\if@bburx \else \no@bbtrue \fi
		\if@bbury \else \no@bbtrue \fi
		\ifno@bb \bb@missing \fi
		\ifno@bb \typeout{FATAL ERROR: no bb supplied or found}
			\no-bb-error
		\fi
		\count203=\@p@sbburx
		\count204=\@p@sbbury
		\advance\count203 by -\@p@sbbllx
		\advance\count204 by -\@p@sbblly
		\edef\@bbw{\number\count203}
		\edef\@bbh{\number\count204}
}
\def\in@hundreds#1#2#3{\count240=#2 \count241=#3
		     \count100=\count240	
		     \divide\count100 by \count241
		     \count101=\count100
		     \multiply\count101 by \count241
		     \advance\count240 by -\count101
		     \multiply\count240 by 10
		     \count101=\count240	
		     \divide\count101 by \count241
		     \count102=\count101
		     \multiply\count102 by \count241
		     \advance\count240 by -\count102
		     \multiply\count240 by 10
		     \count102=\count240	
		     \divide\count102 by \count241
		     \count200=#1\count205=0
		     \count201=\count200
			\multiply\count201 by \count100
		     	\advance\count205 by \count201
		     \count201=\count200
			\divide\count201 by 10
		     	\multiply\count201 by \count101
			\advance\count205 by \count201
		     \count201=\count200
			\divide\count201 by 100
			\multiply\count201 by \count102
			\advance\count205 by \count201
		     \edef\@result{\number\count205}
}
\def\compute@wfromh{
		\in@hundreds{\@p@sheight}{\@bbw}{\@bbh}
		\edef\@p@swidth{\@result}
}
\def\compute@hfromw{
		\in@hundreds{\@p@swidth}{\@bbh}{\@bbw}
		\edef\@p@sheight{\@result}
}
\def\compute@handw{
		\if@height
			\if@width
			\else
				\compute@wfromh
			\fi
		\else
			\if@width
				\compute@hfromw
			\else
				\edef\@p@sheight{\@bbh}
				\edef\@p@swidth{\@bbw}
			\fi
		\fi
}
\def\compute@resv{
		\if@rheight \else \edef\@p@srheight{\@p@sheight} \fi
		\if@rwidth \else \edef\@p@srwidth{\@p@swidth} \fi
		\edef\@p@srheight@total{\@p@srheight}
}
\newtoks\@pshape
\def\@c@ttoplines{\count120}
\def\@c@theightcount{\count121}
\def\@c@tshapecount{\count122}
\newdimen\@d@mwidthshape
\newdimen\@d@mextraheight
\newdimen\@d@mscratch
\def\compute@parshape{
	\if@right
		\if@left
	   		\typeout{error: Can't have both left and right set}
			\@leftfalse
		\fi
	\fi
	\@d@mscratch=\@p@swidth truesp
	\divide \@d@mscratch by 19
	\multiply \@d@mscratch by 20
	\edef\@p@swidthdimen{\the\@d@mscratch}
	\@c@tshapecount=\@c@ttoplines
 	\@d@mscratch=\baselineskip
	\multiply \@d@mscratch by \@c@ttoplines
	\advance \@d@mscratch by .4\baselineskip
    	\edef\@p@stopdistance{\the\@d@mscratch }
	\@d@mscratch=\@p@sheight truesp
	\divide \@d@mscratch by 2
	\edef\@p@shalfboxheight{\the\@d@mscratch}
	\if@toplines
		\loop \@pshape=\expandafter{\the\@pshape 0pt \hsize}
		\advance\@c@ttoplines by -1
		\ifnum\@c@ttoplines>0 \repeat
	\fi
   	\@c@theightcount=\@p@srheight@total
	\advance \@c@theightcount by \@d@mextraheight
	\divide  \@c@theightcount by \baselineskip
	\advance \@c@theightcount by 1
    	\advance \@c@tshapecount by \@c@theightcount
	\advance \@c@theightcount by -1
	\@d@mwidthshape=\hsize
     	\advance \@d@mwidthshape by -\@p@swidthdimen
	\if@left
		\def\@moveshape{0pt}
		\ifnum\@c@theightcount>0
		      	\loop
			\@pshape=%
			\expandafter{\the\@pshape %
					\@p@swidthdimen \@d@mwidthshape}
			\advance \@c@theightcount by -1
			\ifnum\@c@theightcount>0 \repeat
		\else
			\advance \@c@tshapecount by 1
		\fi
	\fi
	\if@right
		\@d@mscratch=\hsize
		\advance \@d@mscratch by -\@p@swidth truesp
		\edef\@moveshape{\@d@mscratch}
		\ifnum\@c@theightcount>0
			\loop
			\@pshape=\expandafter{\the\@pshape 0pt \@d@mwidthshape}
			\advance \@c@theightcount by -1
			\ifnum\@c@theightcount>0 \repeat
		\else
			\advance \@c@tshapecount by 1
		\fi
	\fi
	\ifnum \@p@srheight=0
		\@pshape={}
		\@c@tshapecount = 0
	\else
	 	\@pshape=\expandafter{\the\@pshape 0pt \hsize}
	\fi
}
\def\@p@ssetsurroundboxes{
	\global\parshape=\count122 \the\@pshape		
 	\moveright\@moveshape
	\vbox to 0pt\bgroup\hskip0pt\vskip\@p@stopdistance
}
\newtoks\@captiontok
\newbox\@b@xcaption
\newdimen\@d@mcaptionwidth
\newdimen\@d@mcaptionheight
\newwrite\@w@rcaption
\newtoks\@w@rname
\def\setcaption#1{\@captiontok={#1}}
\def\@set@caption{
	\setbox\@b@xcaption\vbox{\hsize\@d@mcaptionwidth
	\tolerance=9000 \baselineskip=11.4pt
	\noindent\relax\the\@captiontok}
	\if@captionwrite
		\if@captionopen
		\else
			\global\@captionopentrue
			\immediate\openout\@w@rcaption=\the\@w@rname
		\fi
		\immediate\write\@w@rcaption{\the\@captiontok}
		\immediate\write\@w@rcaption{}
	\fi
}
\def\compute@caption{
	\if@captionwidth
	\else
		\@d@mcaptionwidth = \@p@swidth truesp
		\divide \@d@mcaptionwidth by 20
		\multiply \@d@mcaptionwidth by 17
	\fi
	\@set@caption
	\@d@mcaptionheight=\ht\@b@xcaption
	\if@rheight
	\else
		\count100=\@p@srheight
	   	\advance \count100 by \@d@mcaptionheight
	   	\advance \count100 by \bigskipamount
	   	\advance \count100 by \medskipamount
	   	\edef\@p@srheight@total{\number\count100}
	\fi
}
\newif\if@alreadyjtem \@alreadyjtemfalse
\def\newpar{\hfil\vadjust{\vskip\parskip}%
	{\count100=\parskip
	\count101=\baselineskip
	\divide\count101 by 10  \multiply\count101 by 3
	\advance \count100 by \count101
	\divide\count100 by \baselineskip
	\advance\count100 by \prevgraf
	\global\prevgraf=\count100}%
	\break\if@alreadyjtem\else\indent\fi%
}
\let\sav@par=\par
\def\jtem#1{%
    	\if@alreadyjtem\else\bgroup\fi
	\def\par{\sav@par\egroup\sav@par}
	\if@alreadyjtem\else\leftskip\parindent\fi
	\@alreadyjtemtrue
	\noindent\hskip0pt
	\llap{#1\ }\ignorespaces
}
\def\compute@sizes{%
	\compute@bb
	\compute@handw
  	\compute@resv
	\if@caption
		\compute@caption
	\fi
	\if@surround
		\compute@parshape
	\fi
}
\def\@p@sdospecials{
	\ifnum\@p@scost<\@psdraft
	       	\typeout{psfig: including \@p@sfile \space }
	\fi
	\special{ps::[begin] 	\@p@swidth \space \@p@sheight \space
			\@p@sbbllx \space \@p@sbblly \space
			\@p@sbburx \space \@p@sbbury \space
			startTexFig \space }
	\ifnum\@p@scost<\@psdraft
		\if@clip
			\typeout{(clip)}
			\special{ps:: \@p@sbbllx \space \@p@sbblly \space
				\@p@sbburx \space \@p@sbbury \space
			    	doclip \space }
		\fi
	\fi
	\if@box
		\typeout{(box)}
  		\special{ps:: \@p@sbbllx \space \@p@sbblly \space
			\@p@sbburx \space \@p@sbbury \space
		    	dobox \space }
	\fi
	\ifnum\@p@scost<\@psdraft
		\if@prologfile
	    		\special{ps: plotfile \@prologfileval \space }
		\fi
		\special{ps: plotfile \@p@sfile \space }
    		\if@postlogfile
			\special{ps: plotfile \@postlogfileval \space }
		\fi
	\fi
	\special{ps::[end] endTexFig \space }
}
\newif\if@putinvbox

\def\psfig#1{{%
	\ifhmode%
		\vbox\bgroup
		\@putinvboxtrue
	\else
		\@putinvboxfalse
	\fi
       	\ps@init@parms
	\parse@ps@parms{#1}
       	\compute@sizes
	\if@surround
		\psfig@for@surround
	\else
		\psfig@for@regular
	\fi
	\if@putinvbox
       		\egroup
	\fi
}}
\def\psfig@for@surround{%
	\@p@ssetsurroundboxes
	\ifnum\@p@scost<\@psdraft
		\@p@sdospecials
		\vbox to \@p@srheight true sp{\vss}
       	\else
		\if@box
			\@p@sdospecials
		\fi
		\vbox to \@p@srheight true sp{
			\vskip\@p@shalfboxheight
			\hbox to \@p@srwidth true sp{
				\hss
				\ifnum\@p@scost<\@psdraftspecial
					\@p@sfile
				\else
					\@p@sfake
				\fi
      				\hss
			}
		\vss
		}
	\fi
	\if@caption
		\medskip
		\hbox to \@p@srwidth true sp{
			\hss
			\box\@b@xcaption
			\hss
		}
 		\medskip
	\fi
	\vss\egroup
	\vskip-\parskip
}

\def\psfig@for@regular{%
	\if@putinvbox
	\else
		\vskip\parskip
	\fi
	\ifnum\@p@scost<\@psdraft
		\@p@sdospecials
		\vbox to \@p@srheight true sp{%
			\hbox to \@p@srwidth true sp{
			\hfil
			}
		\vfil
		}
       	\else
		\if@box
			\@p@sdospecials
		\fi
	    	\vbox to \@p@srheight true sp{
			\vss
			\hbox to \@p@srwidth true sp{
				\hss
				\ifnum\@p@scost<\@psdraftspecial
					\@p@sfile
				\else
					\@p@sfake
				\fi
				\hss
			}
		    	\vss
		}
	\fi
	\if@caption
		\medskip
		\hbox to \@p@srwidth true sp{
			\hss
			\box\@b@xcaption
			\hss
		}
		\bigskip
	\fi
	\if@putinvbox
	\else
		\vskip-\parskip
	\fi
}
\catcode`\@=12\relax

\psfiginit

\font\scriptsizebbfont=msbm7 scaled \magstep 1
 1
\font\footnotesizebbfont=msbm9 scaled \magstep 0
\font\smallbbfont=msbm7 scaled \magstep 2
\font\bbfont=msbm9 scaled \magstep1  
 1
 2
 3
 4

\def\scriptsizeBbb#1{\hbox{\scriptsizebbfont #1}}

\def\footnotesizeBbb#1{\hbox{\footnotesizebbfont #1}}
\def\smallBbb#1{\hbox{\smallbbfont #1}}
\def\Bbb#1{\hbox{\bbfont #1}}

\newcommand{\Eq}{\mbox{\rm Eq}}

\newcommand{\NE}{\mbox{\it NE}\,}
\newcommand{\PGL}{\mbox{\it PGL}\,}

\newcommand{\Sing}{\mbox{\it Sing}\,}

\newcommand{\Spec}{\mbox{\it Spec}\,}

\newcommand{\SU}{\mbox{\it SU}\,}

\newcommand{\bboxtimes}{\Box\hspace{-1.75ex}\raisebox{.15ex}{$\times$}\,}

\newcommand{\ev}{\mbox{\it ev}\,}

\newcommand{\pt}{\mbox{\it pt}}

\newcommand{\rel}{\mbox{\it\scriptsize rel}}

\newcommand{\virt}{\mbox{\scriptsize\it virt}}

\begin{document}

\enlargethispage{23cm}

\begin{titlepage}

$ $

\vspace{-2cm} 

\noindent\hspace{-1cm}
\parbox{6cm}{\small October 2004}\
   \hspace{6.5cm}\
   \parbox{5cm}{math.AG/0411038}

\vspace{1.5cm}

\centerline{\large\bf
 Extracting Gromov-Witten invariants of a conifold
}
\vspace{1ex}
\centerline{\large\bf
 from semi-stable reduction and relative GW invariants of pairs
}

\vspace{1.5cm}
\centerline{\large
  Chien-Hao Liu
  \hspace{1ex} and \hspace{1ex}
  Shing-Tung Yau
}

\vspace{2em}
\centerline{\small
 ({\it From C.-H.L.\ $:$ In memory of Barbara M.\ Willman.})}
\vspace{2em}

\begin{quotation}
\centerline{\bf Abstract}
\vspace{0.3cm}
\baselineskip 12pt  
{\small
The study of open/closed string duality and large $N$ duality suggests
 a Gromov-Witten theory for conifolds that sits on the border of
 both a closed Gromov-Witten theory and an open Gromov-Witten theory.
The work of Jun Li on Gromov-Witten theory for a projective singular
 variety of the gluing form $Y_1\cup_D Y_2$, where $D$ is a smooth
 divisor on smooth $Y_1$ and $Y_2$, suggests two methods to study
 Gromov-Witten invariants for a projective conifold:
  one by a direct generalization of his construction to the conifold
   singularity and
 the other by an appropriate semi-stable reduction of a degeneration
  to a conifold and then apply his results on this new degeneration
  to extract Gromov-Witten invariants of the original conifold.
In this work we carry out the second method.
Suggested by the semi-stable reduction,
 we associate to a conifold $Y$ with singular locus
 $\{p_1,\,\ldots\,\}$ a set of smooth variety-divisor pairs
   $(\widetilde{Y},E)$, $(Y_i,D_i)$, $i=1,\,\ldots$, and
  a canonical morphism
   $\widetilde{Y}\cup_{\coprod_i D_i}\coprod_i Y_i \rightarrow Y$,
 where
  $\widetilde{Y}$ is the blow-up of $Y$ at the conifold singularities,
  $E\simeq \coprod_i D_i$ is the exceptional divisor, and
  $Y_i$ is a smooth quadric hypersurface in ${\smallBbb P}^4$ with
   a smooth hyperplane section
   $D_i\simeq {\smallBbb P}^1\times{\smallBbb P}^1$, and
  $i$ runs through the labels of the conifold singularities $p_i$.
The existence of a ${\smallBbb Z}/2{\smallBbb Z}$-action on
 the quadric hypersurface $Y_i$ that restricts to an exchange of
 the two product factor ${\smallBbb P}^1$'s of $D_i$ implies that
 the variety $\widetilde{Y}\cup_{\coprod_i D_i}\coprod_i Y_i$
 from gluing these pairs is uniquely determined by $Y$ up to
 isomorphisms.
Jun Li's relative Gromov-Witten theory and invariants for smooth
 variety-divisor pairs $(Y_0, \coprod_iD_i)$ and $(Y_i,D_i)$ and
 a refinement of his degeneration formula worked out here for
 the current situation are then employed to extract Gromov-Witten
 invariants of $Y$.
} 
\end{quotation}

\vspace{2.4em}

\baselineskip 12pt
{\footnotesize
 \noindent
 {\bf Key words:} \parbox[t]{13cm}{
  string world-sheet instanton, open/closed string duality, conifold,
  stable morphism, Gromov-Witten invariant, semi-stable reduction,
  quadric hypersurface, degeneration formula.}
} 

\medskip

\noindent {\small
MSC number 2000$\,$:
 14N35, 81T30.
} 

\bigskip

\baselineskip 11pt
{\footnotesize
\noindent{\bf Acknowledgements.}
 We thank
  Joe Harris, Jun Li, Kefeng Liu, and Cumrun Vafa
 for many influences on this work via lectures, communications,
 and/or very profound works.
 C.-H.L.\ would like to thank in addition
  K.L.\ for serving as a sounding board to many of his raw ideas;
  Lisa Randall, Itay Yavin for valuable lectures;
  Andrew Strominger, Xiaowei Wang for discussions;
  Department of Mathematics of UCLA for hospitality;
  Rev.\ Campbell Willman and Ann Willman for a conversation;
 and Ling-Miao Chou for the tremendous moral support.
 The work is supported by NSF grants DMS-9803347 and DMS-0074329.
} 

\end{titlepage}

\newpage
$ $

\vspace{-4em}  

\centerline{\sc
 Extracting Gromov-Witten Invariants of Conifolds from Pairs
}

\vspace{2em}

\baselineskip 14pt  

\begin{flushleft}
{\Large\bf 0. Introduction and outline.}
\end{flushleft}

\begin{flushleft}
{\bf Introduction:
     Gromov-Witten theory for conifolds in an open/closed string duality
     and a large $N$ duality.}
\end{flushleft}
Given the $3$-sphere $S^3$, let
 $X_0$ be the conifold from the degeneration of $X:=T^{\ast}S^3$,
  as a complex $3$-fold, that pinches the zero-section $3$-cycle
  $S^3$ of $T^{\ast}S^3$ and
 $X^{\prime}$
  be the complex $3$-fold from the small resolution of $X_0$
   with exceptional locus $\simeq {\Bbb P}^1$.
 ($X^{\prime}$ is isomorphic to the total space of the bundle
   ${\cal O}_{{\scriptsizeBbb P}^1}(-1)
                    \oplus{\cal O}_{{\scriptsizeBbb P}^1}(-1)$
  and is a non-compact Calabi-Yau $3$-fold.)
Gopakumar and Vafa [G-V] conjecture the following correspondence
 that relates the $U(N)$ or $\SU(N)$
 Chern-Simons gauge theory on $S^3$ and
 an A-model topological closed string theory on $X^{\prime}$
 via the mechanisms indicated below:

$$
 \begin{array}{ccccc}
  & \parbox{10ex}{\raggedright\scriptsize
     't Hooft expression of Feynman diagrams}
   & & \parbox{15ex}{\scriptsize
      sum over holes/ \newline
      boundaries on the string world-sheet\newline} &  \\[-2ex]
  \framebox[16ex]{$\;$\parbox{15ex}{\footnotesize
      $U(N)$ or $\SU(N)$ \newline
      Chern-Simons \newline gauge theory \newline on $S^3$}}
   & \begin{array}{c}\Longrightarrow\\[-.6ex] \Longleftarrow\end{array}
   & \framebox[19ex]{\hspace{2.4ex}\parbox{20ex}{\footnotesize
      A-model topological \newline open string theory \newline
      on $X$ with boundary \newline on stacked D-brane \newline
      wrapped on $S^3$}}
   & \begin{array}{c}\Longrightarrow\\[-.6ex] \Longleftarrow\end{array}
   &
   \framebox[19ex]{\hspace{4ex}\parbox{21ex}{\footnotesize
      A-model topological \newline closed string theory \newline
      on $X^{\prime}$}} \\[-1ex]
  & \parbox{11ex}{\raggedright\scriptsize
      induced low- \newline energy theory \newline on D-brane world-volume}
   & & \parbox{14ex}{\scriptsize
      $C$-domain field \newline integration \newline to generate holes} &
 \end{array}
$$

Under this correspondence,

$$
\mbox{\footnotesize
\begin{tabular}{lcl}
 \parbox{15em}{$U(N)$ or $SU(N)$ Chern-Simons \newline
               gauge theory on $S^3$}
  & $\Longleftrightarrow$
  & \parbox{15em}{A-model topological closed \newline
                  string theory on $X^{\prime}$}
                              \\[1.6ex] \cline{1-1} \cline {3-3} \\
 $\cdot$ 't Hooft coupling $\lambda := Ng_s$
  & & $\cdot$ $B$-field magnitude on ${\footnotesizeBbb P}^1$ \\[.6ex]
 $\cdot$ Wilson loop observable
  & & $\cdot$
      \parbox[t]{15em}{quantity in the effective theory of\newline
          associated brane-probe in $X^{\prime}$} \\[3ex]
 $\cdot$
 \parbox[t]{15em}{large $N$ limit of Chern-Simons \newline
                  gauge theory on $S^3$}
  & & $\cdot$
      \parbox[t]{15em}{A-model topological string theory \newline
                       on conifold $X_0$}
\end{tabular}
} 
$$

\bigskip

\noindent
This diagram involves
 a large $N$ duality and an open/closed string duality,
and an A-model topological string theory with a conifold as the target
 space serves as a geometric mediator and transition-point for these
 dualities.
The dualities were tested/examined from {\it five} different aspects of
 stringy dualities:
  the {\it string world-sheet} aspect,
  the {\it target space-time} aspect,
  the {\it low dimensional effective field theory} aspect,
  the {\it brane-probe} aspect, and
  the {\it Wilson's theory space} aspect.
See [G-P], [G-V], [O-V1], [O-V2], [Va], [Wi1], [Wi2]
 for string-theoretical details and insights,
 [G-R] for an introductory mathematical review and more references,
 and e.g.\ [A-M-V] and [D-F-G] for generalizations to more general
 non-compact Calabi-Yau $3$-folds involving toric geometry.

Here a {\it conifold} $X_0$ (e.g.\ [B-L], [C-dlO-G-P], and [St]) is
 by definition the singular variety from a degeneration of Calabi-Yau
 $3$-fold $X$ via a deformation of complex structures that pinches
 isolated smoothly-embedded $3$-spheres $S^3$ in a smooth $3$-fold $X$.
And the partition function of an A-model topological string theory with
 target $X_0$ is supposed to compute the string world-sheet instanton
  numbers that can be interpreted either as a counting of holomorphic
  maps from (complex) curves to $X_0$ or from bordered Riemann
  surfaces to $X_0$ with boundary components, if non-empty, mapped to
  the isolated singularities of $X_0$.
In other words, this is a {\it Gromov-Witten theory for conifolds} on
 the mathematical side.
And the open/closed string duality from string theory reveals a very
 distinguished feature of it as a theory on the border of {\it both}
 a closed Gromov-Witten theory on one geometry and
 an {\it open} Gromov-Witten theory on another geometry related to
  the previous geometry via the conifold transition.

To develop a Gromov-Witten theory for conifolds, there are two paths
 one may attempt to follow.
The first one is to generalize the techniques in [Li1] and [Li2]
 directly to a conifold singularity. This is technical.
The second one is to try to replace a conifold as the degenerate
 fiber of a smooth ${\Bbb A}^1$-family $W/{\Bbb A}^1$ by the
 degenerate fiber of a semi-stable reduction of $W/{\Bbb A}^1$ and
 see if one can reduce the problem to the case already dealt with in
 [Li1] and [Li2] and use it to extract the Gromov-Witten invariants
 of the conifold.
Surely, for general singularities one will not expect the second path
 would immediately work either since one still misses the understanding
 of Gromov-Witten theory for a pair $(Y,D)$, where $Y$ is smooth and $D$
 is a divisor on $Y$ with simple normal crossing singularities,
 and a related degeneration/gluing formula for gluing Gromov-Witten
 invariants from a collection of such pairs.
As if given by God, it turns out that such potentially existing
 difficulties along the second path are not really there
 for a conifold singularity, (Sec.~1).
This accidental simplicity for conifolds together with similar
 discussions and results in \mbox{[L-Y]} that refine Jun Li's
 degeneration formula to one for curve classes in
 $H_2(\,\cdot\,;{\Bbb Z})$ or $A_1(\,\cdot\,)$, (Sec.~2),
 enable us to extract some Gromov-Witten invariants of conifolds,
 (Sec.~3).
We explain the details of this second path in this work.

\bigskip

\noindent
{\it Convention.}
 This work is parallel to [L-Y] and follows the notations and
  the terminology of [Li1], [Li2], and [L-Y] closely, except
  where confusions may occur.
 Other notations follow [Hart], [Fu], [De], and [K-M].
 All schemes are over ${\Bbb C}$ and all points are referred
  to closed points.
 All conifolds are assumed to be projective.
 Though appearing a few times in the Introduction for easy match
  with literatures in string theory, the Calabi-Yau condition is
  not relevant in this work but will be relevant in an application.

\bigskip

\noindent
{\it Preliminary.}
 Readers are referred to [Li1: Sec.~0], [Li2: Sec.~0], and
  [L-Y: Sec.~1, Sec.~2] for definitions and an outline of Jun Li's
  work that are needed for the current work.

\bigskip

\begin{flushleft}
{\bf Outline.}
\end{flushleft}
{\small
\baselineskip 11pt  
\begin{quote}
%
 1. A semi-stable reduction of a conifold degeneration.

 2. A degeneration formula of Gromov-Witten invariants
    with respect to a curve class.

 3. Extracting Gromov-Witten invariants of a conifold from pairs.
\end{quote}
} 

\bigskip

\baselineskip 14pt  


\newpage
\section{A semi-stable reduction of a conifold degeneration.}

\begin{flushleft}
{\bf Semi-stable reduction and the associated smooth pairs of
     a conifold.}
\end{flushleft}
Let $\pi:W\rightarrow {\Bbb A}^1$ be a flat family of (complex)
 $3$-dimensional projective varieties with smooth general fibers $W_t$,
  $t\ne {\mathbf 0}$, and smooth total space $W$ such that 
 $W_t$ degenerates to a conifold $W_0=Y$ over ${\mathbf 0}\in{\Bbb A}^1$.
For simplicity of notations and presentations, we assume that $Y$ has
 only one conifold singularity.
An example of such family can be obtained from a degeneration of
 quintic $3$-folds in ${\Bbb P}^4$, cf.\ [C-dlO-G-P].
In a local analytic germ, the degeneration to a conifold singularity
 is modeled on the morphism
 $$
  \begin{array}{ccc}
   \Spec{\Bbb C}[[\,x,y,z,w\,]]
    & \longrightarrow   & \Spec {\Bbb C}[[t]] \\[.6ex]
   xy-zw                & \longleftarrow    & t \,.
  \end{array}
 $$
This gives a $3$-dimensional isolated hypersurface singularity of
 multiplicity $2$ at the fiber over ${\mathbf 0}:=(t)$.
The family $W/{\Bbb A}^1$ can be semi-stabilized by the following
 sequence of blow-ups and finite base change obtained
 from a straightforward computation:
 $$
  \begin{array}{cccccccl}
   W^{(3)}  & \stackrel{\varphi_3}{\longrightarrow}  & W^{(2)}
            & \stackrel{\varphi_2}{\longrightarrow}  & W^{(1)}
            & \stackrel{\varphi_1}{\longrightarrow}  & W^{(0)} := W \\
   \mbox{\scriptsize $\pi_3$}\downarrow\hspace{3ex}
    && \mbox{\scriptsize $\pi_2$}\downarrow\hspace{3ex}
    && \hspace{2ex}\downarrow\mbox{\scriptsize $\pi_1$}
    && \hspace{6ex}\downarrow\mbox{\scriptsize $\pi_0=\pi$}\\[.6ex]
   {\Bbb A}^1   & =  & {\Bbb A}^1   & \stackrel{\alpha}{\longrightarrow}
    & {\Bbb A}^1     & =  & {\Bbb A}^1   &,
  \end{array}
 $$
 where
 \begin{itemize}
  \item[$\cdot$]
   $\varphi_1: W^{(1)}\rightarrow W^{(0)}$ is the blow-up of $W$
   at the conifold singularity of $W_0\,$;

  \item[$\cdot$]
   recall that $W$ is smooth and hence the exceptional locus of
    $\varphi_1$ is a ${\Bbb P}^3\,$;
   the degenerate fiber $\pi_1^{-1}({\mathbf 0})$ of $\pi_1$
    contains a nonreduced irreducible component of multiplicity $2$
    supported on this ${\Bbb P}^3\,$;
   as divisors in $W^{(1)}\,$,
   $\pi_1^{-1}({\mathbf 0})=\widetilde{Y}+ 2\,{\Bbb P}^3$,
   where
    $\widetilde{Y}$ is the resolution of $Y$ by a blow-up at the
     conifold singularity and
    $\widetilde{Y}\cap {\Bbb P}^3$ is a smooth quadric surface in
     ${\Bbb P}^3\,$;

  \item[$\cdot$]
   $\alpha:({\Bbb A}^1,{\mathbf 0})\rightarrow({\Bbb A}^1,{\mathbf 0})$
   is a finite morphism of degree $2$ branched over ${\mathbf 0}$,
   $\varphi_2$ and $\pi_2$ on $W^{(2)}$ are from the fibered-product
    of $\alpha$ and $\pi_1$;

  \item[$\cdot$]
   $W^{(2)}$ is now a singular scheme whose singularities
   are modelled on the Whitney's umbrella;
   the singular locus $\Sing(W^{(2)})$ of $W^{(2)}$ with the reduced
    subscheme structure has multiplicity $2$, lies over ${\mathbf 0}$,
    and is isomorphic to $\Bbb P^3\,$;

  \item[$\cdot$]
   normalization/blow-up of $W^{(2)}$ along
    $\Sing(W^{(2)})\simeq{\Bbb P}^3$ gives $W^{(3)}$ which is smooth
    with $\pi_3^{-1}({\mathbf 0})=Y_0\cup Y_1$,
    where $Y_0:= \widetilde{Y}$ and $Y_1$ is naturally realized as
     a double cover of ${\Bbb P}^3$ branched over a smooth quadric
     surface in ${\Bbb P}^3$.
 \end{itemize}

In local germs or local formal schemes,
 the morphisms in the diagram above are given by

{\footnotesize
 $$
  \begin{array}{cccl}
   \Spec{\footnotesizeBbb C}[[x_1,y_1,z_1,w_1]]\hspace{5em}
     & \stackrel{\varphi_1}{\longrightarrow}
     & \Spec{\footnotesizeBbb C}[[x_0,y_0,z_0,w_0]]  \\
   w_1^2(x_1y_1-z_1)
    \hspace{4em} x_1w_1\,;\; y_1w_1\,;\; z_1w_1\,;\; w_1
    & \stackrel{\varphi_1^{\sharp}}{\longleftarrow}
    & x_0\,;\; y_0\,;\; z_0\,;\; w_0 \hspace{4.6em} x_0y_0-z_0w_0
                                                \hspace{1.4em} \\[.6ex]
   \mbox{\scriptsize $\pi_1^{\sharp}$}\uparrow
    \hspace{3.2em}\mbox{\scriptsize $\pi_1$}\downarrow\hspace{11em}
   & & \hspace{6em}\downarrow\mbox{\scriptsize $\pi_0$}
        \hspace{3em}\uparrow\mbox{\scriptsize $\pi_0^{\sharp}$}  \\[.6ex]
   t_1 \hspace{12em} t_1\hspace{2em}
    & \longleftarrow  & \hspace{1ex}t_0 \hspace{9em}  t_0     \\
   \Spec{\footnotesizeBbb C}[[t_1]] \hspace{4em}
    & \longrightarrow & \Spec{\footnotesizeBbb C}[[t_0]] &,
  \end{array}
 $$

 $$
  \begin{array}{cccl}
   \Spec{\footnotesizeBbb C}[[x_1,y_1,z_1,w_1, t_2]]
                                  /(t_2^2-w_1^2(x_1y_1-z_1))
     & \stackrel{\varphi_2}{\longrightarrow}
     & \Spec{\footnotesizeBbb C}[[x_1,y_1,z_1,w_1]]  \\
   t_2
    \hspace{8em} x_1\,;\; y_1\,;\; z_1\,;\; w_1
    & \stackrel{\varphi_2^{\sharp}}{\longleftarrow}
    & x_1\,;\; y_1\,;\; z_1\,;\; w_1 \hspace{4.6em} w_1^2(x_1y_1-z_1)
                                                \hspace{1.4em} \\[.6ex]
   \mbox{\scriptsize $\pi_2^{\sharp}$}\uparrow
    \hspace{3.2em}\mbox{\scriptsize $\pi_2$}\downarrow\hspace{11em}
   & & \hspace{6em}\downarrow\mbox{\scriptsize $\pi_1$}
        \hspace{3em}\uparrow\mbox{\scriptsize $\pi_1^{\sharp}$}  \\[.6ex]
   t_2 \hspace{12em} t_2^2\hspace{2em}
    & \longleftarrow  & \hspace{1ex}t_1 \hspace{9em}  t_1     \\
   \Spec{\footnotesizeBbb C}[[t_2]] \hspace{4em}
    & \longrightarrow & \Spec{\footnotesizeBbb C}[[t_1]] &;
  \end{array}
 $$
}

\noindent
consider the ring isomorphism generated by
 $x_1\mapsto x_2$, $y_1\mapsto y_2$, $z_1\mapsto x_2y_2-z_2$,
   $w_1\mapsto w_2$, and $t_2\mapsto t_2$
and rewrite
$$
 \mbox{\footnotesize
  $\Spec{\footnotesizeBbb C}[[x_1,y_1,z_1,w_1,t_2]]/
                           (t_2^2-w_1^2(x_1y_1-z_1))$}
 \hspace{2em}\mbox{as}\hspace{2em}
 \mbox{\footnotesize
  $\Spec{\footnotesizeBbb C}[[x_2,y_2,z_2,w_2,t_2]]/
                           (t_2^2-w_2^2z_2)$}\,,
$$
 which reveals the Whitney umbrella transverse surface singularity
  of $W^{(2)}$;
 and then

{\footnotesize
 $$
  \begin{array}{cccl}
   \Spec{\footnotesizeBbb C}[[x_3,y_3,z_3,w_3]] \hspace{6em}
     & \stackrel{\varphi_3}{\longrightarrow}
     & \Spec{\footnotesizeBbb C}[[x_2,y_2,z_2,w_2,t_2]]/
                                               (t_2^2-w_2^2z_2)  \\
   z_3w_3 \hspace{6em} x_3\,;\; y_3\,;\; z_3^2\,;\; w_3\,; z_3w_3
    & \stackrel{\varphi_3^{\sharp}}{\longleftarrow}
    & x_2\,;\; y_2\,;\; z_2\,;\; w_2\,;t_2 \hspace{10em} t_2
                                                \hspace{2em} \\[.6ex]
   \mbox{\scriptsize $\pi_3^{\sharp}$}\uparrow
    \hspace{3.2em}\mbox{\scriptsize $\pi_3$}\downarrow\hspace{12em}
   & & \hspace{11em}\downarrow\mbox{\scriptsize $\pi_2$}
        \hspace{4em}\uparrow\mbox{\scriptsize $\pi_2^{\sharp}$}  \\[.6ex]
   t_3 \hspace{12em} t_3\hspace{2.6em}
    & \longleftarrow  & \hspace{3em}t_2 \hspace{12.6em}  t_2     \\
   \Spec{\footnotesizeBbb C}[[t_3]] \hspace{5em}
    & \longrightarrow & \hspace{4.4em}\Spec{\footnotesizeBbb C}[[t_2]] &.
  \end{array}
 $$
} 

\noindent
(Cf.\ {\sc Figure}~1-1.)
%

\begin{figure}[htbp]
 \setcaption{{\sc Figure} 1-1.
  \baselineskip 14pt
  A semi-stable reduction of a conifold degeneration.
  The degenerate fiber over ${\mathbf 0}\in{\Bbb A}^1$
   in each family is indicated.
 } 
\centerline{\psfig{figure=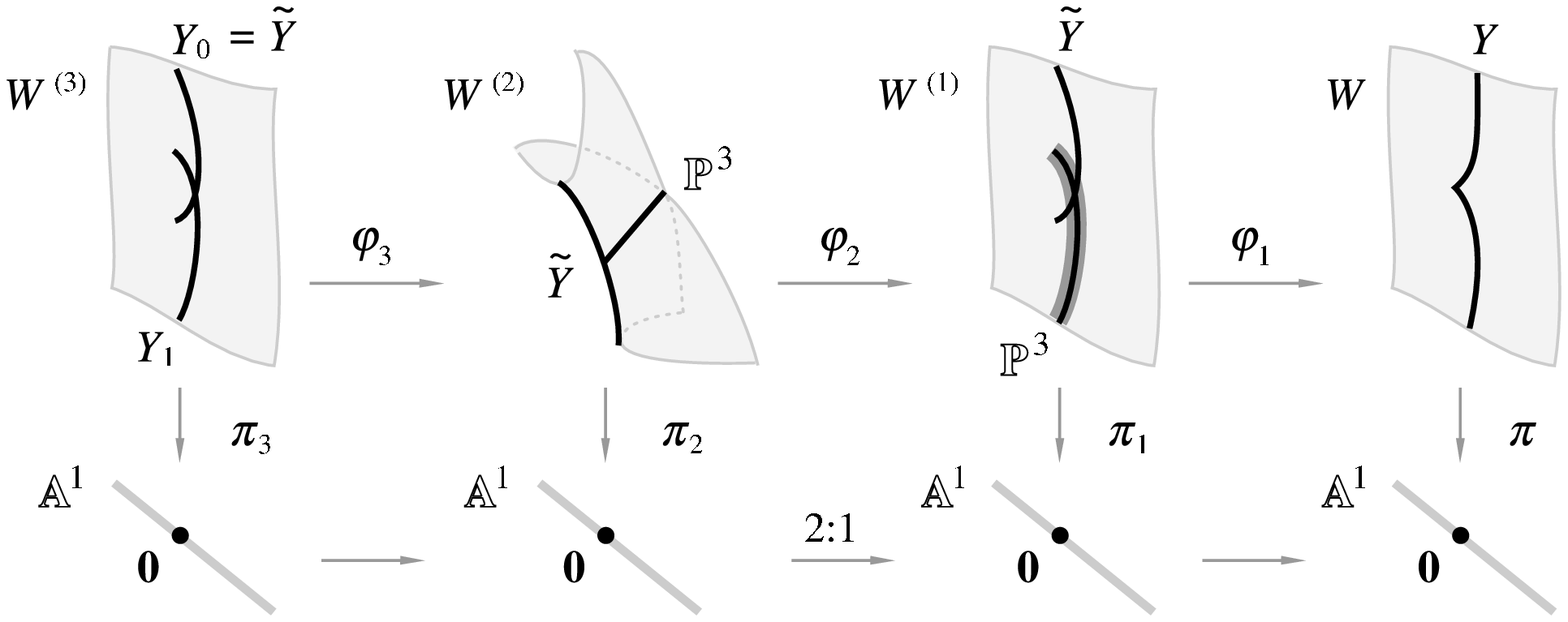,width=13cm,caption=}}
\end{figure}

\bigskip

\noindent
{\bf Definition 1.1 [canonical semi-stable reduction].} {\rm
 We shall call $W^{(3)}/{\Bbb A}^1$ obtained above
  the {\it canonical semi-stable reduction} of the conifold
  degeneration $W/{\Bbb A}^1$.
} 

\bigskip

After this semi-stable reduction of the conifold degeneration
 $\pi:W \rightarrow {\Bbb A}^1$, one obtains a new family
  $\pi_3:W^{(3)}\rightarrow {\Bbb A}^1$ of $3$-folds with
   $W^{(3)}$ smooth,
   $W^{(3)}_t=W_{t^2}$ for $t\ne {\mathbf 0}$, and
   $W^{(3)}_0 := \pi_3^{-1}({\mathbf 0})=Y_0\cup_D Y_1$
    is a variety with simple normal crossing singularity from
    gluing smooth $Y_0$ and $Y_1$ along the isomorphic smooth
    divisor $D$.
This happens to be exactly the type of degenerations whose
 Gromov-Witten theory are studied in [Li1] and [Li2].

\bigskip

\noindent
{\bf Definition 1.2 [associated canonical pairs to conifold].} {\rm
 We shall call the set of pairs, $(Y_0,D)$ and $(Y_1,D)$, the
 {\it set of canonical smooth pairs associated to the conifold} $Y$.
} 

\bigskip

\noindent
By construction, there is a canonical morphism
 $p=p_0\cup_D p_1: Y_0\cup _D Y_1\rightarrow Y$,
 where
  $p_0: Y_0\rightarrow Y$ is the resolution of $Y$
   whose exception locus is the smooth divisor $D$ and
  $p_1:Y_1\rightarrow Y$ pinches $Y_1$ to the conifold singularity
   of $Y$.

\bigskip

\begin{flushleft}
{\bf Special properties of $Y_1$.}
\end{flushleft}
The $3$-fold $Y_1$ is naturally realized as a quadric hypersurface
 in ${\Bbb P}^4$ with $D\simeq {\Bbb P}^1\times{\Bbb P}^1$ realized
 as the intersection of $Y_1$ with a hyperplane $H$ in ${\Bbb P}^4$.
It follows from the Bruhat cellular decomposition of a smooth quadric
 hypersurface in a projective space and the Lefschetz hyperplane
 theorem that the only non-vanishing Chow groups and (singular)
 homology groups for $Y_1$ coincide and are given by
 $$
  A_k(Y_1)\; \simeq\; H_{2k}(Y_1;{\Bbb Z})\; \simeq\; {\Bbb Z}\,,
  \hspace{1em}\mbox{for $k=0,\,1,\,2,\,3$}\,.
 $$
In particular,
 $A_1(Y_1)\simeq H_2(Y_1;{\Bbb Z})$ is generated by a complex line
  $\gamma_2$ from either of the rulings of $D$,
 $A_2(Y_1)\simeq H_4(Y_1;{\Bbb Z})$ is generated by the complex surface
  $\gamma_4:=D$, and $\gamma_2\cap\gamma_4=+1$ on $Y_1$.
(See [Fu], [G-H], [Harr], and [P-S] for details.)

\bigskip

Fix
 an embedding $D=Y_1\cap H\subset Y_1\hookrightarrow {\Bbb P^4}$ and
 an isomorphism $D\simeq {\Bbb P}^1\times{\Bbb P}^1$ from the rulings
  of $D$.
Up to a $\PGL(5,{\Bbb C})$-action on ${\Bbb P}^4$, one can write
 this explicitly in terms of homogeneous coordinates
 $[u:x:y:z:w]$ of ${\Bbb P}^4$ as
$$
 D\,=\,V(u\,,\, u^2+xy-zw)\;\subset\; Y_1\,=\,V(u^2+xy-zw)\;
   \subset\; {\Bbb P}^4
 \hspace{1em}\mbox{with}\hspace{1em}H=V(u)\,.
$$
The ${\Bbb Z}/2{\Bbb Z}$-action on $D$ that interchanges the two
 ${\Bbb P}^1$-factors extends to a linear ${\Bbb Z}/2{\Bbb Z}$-action
 on ${\Bbb P}^4$ that leaves $Y_1$ invariant.
Explicitly,
 $D\simeq {\Bbb P}^1\times{\Bbb P}^1$ via
  $\{\,xy-zw=0\,\}\rightarrow ([x:z],[x:w])$
 and one can choose this ${\Bbb Z}/2{\Bbb Z}$-action to be
 the exchange $z\leftrightarrow w$.
It follows that:

\bigskip

\noindent
{\bf Lemma 1.3 [uniqueness of gluing].} {\it
 Any variety from gluing the smooth pairs, $(Y_0,D)$ and $(Y_1,D)$,
 associated to $Y$ via an automorphism of $D$ is isomorphic to
 each other.
} 

\bigskip

{\samepage
\begin{flushleft}
{\bf The monodromy of $\pi:W\rightarrow {\Bbb A}^1$.}
\end{flushleft}
To be able to refine Jun Li's degeneration formula from one with}
 respect to a $\pi_3$-ample line bundle on $W^{(3)}/{\Bbb A}^1$
 to one with respect to a curve class of a general fiber $W_t$ of
 $W/{\Bbb A}^1$, the monodromy around the degenerate fiber
 $W_0=Y$ of $W/{\Bbb A}^1$ is required to be well-behaved.
Let us thus take a look at this.
In this part, we will leave the algebro-geometric category over
 ${\Bbb C}$ and enter the topological category over ${\Bbb R}$.
In particular, `$\simeq$' means `topologically homeomorphic'.

Topologically, the degeneration $W/{\Bbb A}^1$ pinches a real
 $3$-cycle realized as a smoothly embedded $3$-sphere $S^3$ in $W_t$,
 $t\ne {\mathbf 0}$, to the conifold singularity of $W_0=Y$.
Over a loop $S^1\hookrightarrow {\Bbb A}^1$ around
 ${\mathbf 0}\in {\Bbb A}^1$
  (here we will think of ${\Bbb A}^1$ as the usual complex plane
   ${\Bbb C}$ with the analytic topology and $S^1$ is an oriented
   circle therein),
 $\pi^{-1}(S^1)$ is the mapping torus associated to a smooth
 automorphism $\phi:W_t\stackrel{\sim}{\rightarrow} W_t$.
$\phi$ can be homotoped to an automorphism, still denoted by $\phi$,
 on $W_t$ that is not the identity map only in a small tubular
 neighborhood $N_{\varepsilon}(S^3/W_t)\simeq S^3\times B^3$
 of $S^3$ in $W_t$, where $B^3$ is the unit (real) $3$-ball and
 there is an orientation on $(S^3, B^3)$ that is compatible with
 the complex structure on $N_{\varepsilon}(S^3/W_t)$.
Indeed, $\phi$ is a {\it generalized Dehn twist} along $S^3$ and,
 up to homotopy, can be described explicitly as follows.
Recall that $S^3$ is the underlying topology of the Lie group
 $\SU(2)$.
Let
 $e$ be its identity element and
 fix a degree $1$ map
  $\phi_e:(B^3,\partial B^{3\,})\rightarrow (S^3,e)$
  that is a homeomorphism from $B^3-\partial B^3$ to $S^3-\{e\}$.
Define
 $\phi:N_{\varepsilon}(S^3/W_t)\rightarrow N_{\varepsilon}(S^3/W_t)$
 by
 $$
  \begin{array}{cccccl}
   \phi\;: & S^3\times B^3
           & \longrightarrow  & S^3\times B^3 \\[.6ex]
    & (g,x)   & \longmapsto   & ( g\cdot\phi_e(x) ,\,x)  &,
  \end{array} 
 $$
 where the $\cdot$ in the definition is the group multiplication of
 $\SU(2)$.
By construction, $\phi$ restricts to the identity map on
 $\partial(N_{\varepsilon}(S^3/W_t))$ and extends by the identity map
 to the automorphism $\phi$ on $W_t$.

Since $H_i(S^3\times B^3,\partial(S^3\times B^3);{\Bbb Z})=0$
 for $i=0,\,1,\,2$,
the monodromy of $\phi$ on $H_i(W_t;{\Bbb Z})$ is the identity map
 for $i=0,\,1,\,2$.
(The monodromy on $H_3(W_t;{\Bbb Z})$ is the much studied 
 Picard-Lefschetz operation
 $\gamma \mapsto \gamma + \langle\gamma, [S^3]\rangle\,[S^3]$,
  where $\langle\,,\,\rangle$ is the intersection pairing on
  $H_3(W_t;{\Bbb Z})$. But we do not need this for this work.)
Since the conifold $W_0=Y$ is topologically the quotient space of
 $W_t$, $t\ne{\mathbf 0}$, by identifying the vanishing $3$-cycle $S^3$
 to a point, $Y$ is homotopically equivalent to
 $W_t\cup_{S^3\simeq \partial B^4} B^4$, where $t\ne {\mathbf 0}$ and
 the $4$-ball $B^4$ is attached to $W_t$ along an isomorphism of
 $\partial B^4$ to the vanishing $3$-cycle $S^3$.
This implies that $H_i(W_0;{\Bbb Z})\simeq\; H_i(W_t;{\Bbb Z})$
 canonically for $i=0,\,1,\,2$.
In particular, $R^{\,\bullet}\pi_{\ast}{\Bbb Z}_W$ has a summand
 that is the constant sheaf on ${\Bbb A}^1$ associated to the group
 $H_2(W_t; {\Bbb Z})=H_2(Y;{\Bbb Z})$.

To summarize,

\bigskip

\noindent
{\bf Lemma 1.4
     [canonical identification and trivial monodromy on $H_2$].}
{\it
 The monodromy of $W/{\Bbb A}^1$ around the conifold $W_0=Y$
  is trivial on $H_2(W_t;{\Bbb Z})$.
 Each $H_2(W_t;{\Bbb Z})$, $t\ne{\mathbf 0}$, is canonically
  isomorphic to $H_2(W_0;{\Bbb Z})=H_2(Y;{\Bbb Z})$.
} 

\bigskip

We now resume to proceed in the domain of algebraic geometry.

\bigskip

\bigskip

\section{A degeneration formula of Gromov-Witten invariants
         with respect to a curve class.}

Once the semi-stable degeneration $W^{(3)}/{\Bbb A}^1$
 of the conifold degeneration $W/{\Bbb A}^1$ is understood,
the second ingredient toward a Gromov-Witten theory for $Y$ via
 extracting information from gluing pairs $(Y_0,D)$ and $(Y_1,D)$
 is the degeneration/gluing formula of Gromov-Witten invariants
 for the family $W^{(3)}/{\Bbb A}^1$ that is intrinsic to curve
 classes on $Y$ (note: {\it not}$\,$ on $W^{(3)}_0$).
The goal of this section is to derive such a formula from [Li2].
The reasoning and the discussion of this section are parallel to
 [L-Y: Sec.~2 and Sec.~3].
Readers are referred to [L-Y] for more explanations in a simpler
 situation and [Li1] and [Li2] for the thorough related details.

\bigskip

\begin{flushleft}
{\bf A refinement of Jun Li's degeneration formula.}
\end{flushleft}
Recall first the definitions of the following:
{\it admissible weighted graph}
 ([Li1: Definition 4.6], [L-Y: Definition 1.1];
{\it admissible triple}
 ([Li1: Definition 4.11], [L-Y: Definition 1.2]),
the integer triple $|\eta|$ for an admissible triple
 $\eta=(\Gamma_0,\Gamma_1,I)$, and
Jun Li's {\it degeneration formula}
 ([Li2: Sec.~0, Theorem 3.15, and Corollary 3.16]).

Fix a $\pi_3$-ample line bundle ${\cal L}$ on $W^{(3)}/{\Bbb A}^1$
 and an ${\cal L}$-degree $d$, then the stack of stable morphisms
 from prestable curves of genus $g$, $n$ marked points, into fibers
 of the universal family of expanded degenerations associated to
 $W^{(3)}/{\Bbb A}^1$ of ${\cal L}$-degree $d$ is a Deligne-Mumford
 stack ${\frak M}({\frak W}^{(3)}, (g,n;d))$ with a perfect
 obstruction theory.
By construction, ${\frak M}({\frak W}^{(3)}, (g,n;d))$ fibers
 naturally over ${\Bbb A}^1$ with the fiber over ${\mathbf 0}$
 denoted by ${\frak M}({\frak W}^{(3)}_0, (g,n;d))$.
The perfect obstruction theory on
 ${\frak M}({\mathfrak W}^{(3)},(g,n;d))$ restricts to
 a perfect obstruction theory on
 ${\mathfrak M}({\mathfrak W}^{(3)}_0,(g,n;d))$.
The Gromov-Witten theory and invariants of the singular variety
 $Y_0\cup_DY_1$ from gluing is defined via the virtual fundamental
 class $[{\mathfrak M}({\mathfrak W}^{(3)}_0,(g,n;d))]^{\virt}$
 on ${\mathfrak M}({\mathfrak W}^{(3)}_0,(g,n;d))$.
On the other hand, for each of the smooth variety-divisor pairs,
 $(Y_i,D)$, Jun Li constructed a general theory of relative
 Gromov-Witten invariants, which consists of
  the moduli stack ${\mathfrak M}({\mathfrak Y}_i^{\rel},\Gamma_i)$
   of stable morphisms of topological type $\Gamma_i$ from prestable
   curves to the fibers of the universal family of the stack
   ${\mathfrak Y}_i^{\rel}$ of expanded relative pairs associated to
   $(Y_i,D)$ and
  a perfect obstruction theory thereon.
For each admissible triple $\eta$ with $|\eta|=(g,n;d)$,
 there is a morphism,
 (the $r$ below is the total root weight of $\Gamma_0$ in $\eta$),
 $$
  \Phi_{\eta}\;:\;
   {\mathfrak M}({\mathfrak Y}_0^{\rel},\Gamma_0)\times_{D^r}\,
   {\mathfrak M}({\mathfrak Y}_1^{\rel},\Gamma_1)\;
   \longrightarrow\; {\mathfrak M}({\mathfrak W}^{(3)}_0,(g,n;d))
 $$
 that is finite \'{e}tale to its image
 ${\mathfrak M}
     ({\mathfrak Y}_0^{\rel}\sqcup{\mathfrak Y}_1^{\rel},\eta)$
 and has degree a combinatorial factor $|\Eq(\eta)|$ from $\eta$,
 cf.\ [Li1: Sec.~4].
Up to the difference from a nonreduced structure that has to be
 taken care ([Li2: Sec.~4.4]),
 ${\mathfrak M}({\mathfrak W}^{(3)}_0,(g,n;d))$
 is the union of the substacks
 \mbox{${\mathfrak M}
       ({\mathfrak Y}_0^{\rel}\sqcup{\mathfrak Y}_1^{\rel},\eta)$}
 in ${\mathfrak M}({\mathfrak W}^{(3)}_0,(g,n;d))$.
Going from ${\mathfrak M}({\mathfrak W}^{(3)}_t,(g,n;d))$,
 $t\ne{\mathbf 0}$, which gives the usual Gromov-Witten theory
 of the smooth projective $W_t$ to
 ${\mathfrak M}({\mathfrak W}^{(3)}_0,(g,n;d))$, which gives
 the Gromov-Witten theory on $Y_0\cup_D Y_1$, and
 recast it to a form from gluing relative Gromov-Witten theory of
  pairs, $(Y_0,D)$ and $(Y_1,D)$,
 gives a degeneration formula that relates Gromov-Witten invariants
 of $W_t$ to a combination of relative Gromov-Witten invariants
 of $(Y_0,D)$ and $(Y_1,D)$ in the degenerate fiber $W_0^{(3)}$.
See [Li1] and [Li2] for the complete technical details.

Since the monodromy on $H_2(W_t;{\Bbb Z})$ around
 ${\mathbf 0}\in {\Bbb A}^1$ is trivial, so is the monodromy
 on $H_2(W^{(3)}_t;{\Bbb Z})$ around ${\mathbf 0}$.
Recall also the canonical isomorphism
 $H_2(W_t;{\Bbb Z})\simeq H_2(Y;{\Bbb Z})$.
Thus, a $\beta^{\,\prime}\in H_2(Y;{\Bbb Z})$ determines
 a unique class in $H_2(W_t;{\Bbb Z})$, $t\ne{\mathbf 0}$, and
 hence a unique class, still denoted by $\beta^{\,\prime}$,
 in $H_2(W^{(3)}_t;{\Bbb Z})$, $t\ne{\mathbf 0}$.
Define the ${\cal L}$-degree ${\cal L}\cdot\beta^{\,\prime}$ of
 $\beta^{\,\prime}\in H_2(Y;{\Bbb Z})$ via this identification.
Let
 $d$ be the ${\cal L}$-degree for a given curve class
  $\beta\in H_2(Y;{\Bbb Z})$ and
 define
  $C_{({\cal L},d)}(Y)
   :=\{\,\beta^{\,\prime}\in H_2(Y;{\Bbb Z})\,:\,
                    {\cal L}\cdot \beta^{\,\prime}=d\,\}$.
Then
 $$
  {\frak M}({\frak W}^{(3)},(g,n;d))\;
   =\; \coprod_{\beta^{\prime}\in C_{({\cal L},\,d)}(Y)}
       {\mathfrak M}({\mathfrak W}^{(3)}, (g,n;\beta^{\,\prime}))\,,
 $$
 where ${\mathfrak M}({\mathfrak W}^{(3)},(g,n;\beta^{\,\prime}))$
  is the stack of stable morphisms of ${\cal L}$-degree $d$ from
  prestable curves of genus $g$ with $n$ marked points to the fibers
  of the universal family of the stack ${\mathfrak W}^{(3)}$ of
  expanded degenerations associated to $W^{(3)}/{\Bbb A}^1$ such
  that after the post-composition with the morphisms
  $$
   {\mathfrak W}^{(3)}/{\Bbb A}^1\;
    \longrightarrow\; W^{(3)}/{\Bbb A}^1\;
    \stackrel{\varphi_1\circ\varphi_2\circ\varphi_3}{\longrightarrow}\;
    W/{\Bbb A}^1\,,
  $$
  the images of the stable morphisms lie in the curve class
  $\beta^{\,\prime}\in H_2(Y;{\Bbb Z})$.
From its definition,
 ${\mathfrak M}({\mathfrak W}^{(3)},(g,n;\beta^{\,\prime}))$
 depends not just on $(g,n,\beta^{\,\prime})$ but also on ${\cal L}$.
(Indeed, since $Y_1$ is mapped to the conifold singularity of $Y$,
 the choice of curve classes from $H_2(Y_0\cup_D Y_1;{\Bbb Z})$
 does depend on ${\cal L}$ in general.)

Recall
 the ${\cal L}$-dependent set $\Omega_{(g,k;d)}$ of admissible triples
  $\eta$ such that $|\eta|=(g,k;d)$ and
 the quotient set $\overline{\Omega}_{(g,k;d)}$
  defined in [Li2: Sec.~0].
Recall also the morphism \linebreak
 $p=p_0\cup_D p_1: Y_0\cup_D Y_1\rightarrow Y$.
For an admissible weighted graph $\Gamma$ for a relative pair,
 let $b(\Gamma):=\sum_{v\in V(\Gamma)} b(v)$. 
Define the {\it $\beta$-compatible subset} of $\Omega_{(g,n;d)}$ by
 $$
  \Omega_{(g,k;\beta)}^{\cal L}\;
   :=\; \left\{\,
    \eta=(\Gamma_0,\Gamma_1,I)\in\Omega_{(g,k;d)}\,\left|\,
       p_{0\ast}b(\Gamma_0)=\beta \right.
      \,\right\}\,.
 $$
Same discussions as in [L-Y: Sec.~2] imply the existence of
 a perfect obstruction theory on the moduli stack
  ${\mathfrak M}({\mathfrak W}^{(3)},(g,n;\beta))$
   over ${\Bbb A}^1$ and
  its fiber ${\mathfrak M}({\mathfrak W}^{(3)}_0,(g,n;\beta))$
   over ${\mathbf 0}$, inherited from those constructed in [Li2].
This gives a well-defined Gromov-Witten theory and Gromov-Witten
 invariants on $Y_0\cup_D Y_1$ associated to
 $\beta\in H_2(Y;{\Bbb Z})$,
Jun Li's degeneration formula [Li2]
 implies then the following degeneration formula:

\bigskip

\noindent
{\bf Lemma 2.1 [Jun Li's degeneration formula].}
(Cf.\ [L-Y: Corollary 2.2] and explanations of
      notations there and in [Li1] and [Li2].)
{\it
 Let
  $\alpha\in H^0_c(R^{\ast}\pi_{3\ast}{\Bbb Q}_{W^{(3)}})^{\times n}$,
   whose restriction to $W^{(3)}_t$ will be denoted by $\alpha(t)$, and
  $\zeta\in A_{\ast}({\mathfrak M}_{g,n})$.
 Denote by $\Psi_{(g,n;\beta)}^{W^{(3)}_t}(\alpha(t),\zeta)$
  the usual Gromov-Witten invariant of $W^{(3)}_t$ associated to
  these data.
 For $\eta\in\overline{\Omega}_{(g,n;\beta)}^{\cal L}$, assume that
  $G_{\eta}^{\ast}(\zeta)
    =\sum_{j\in K_{\eta}}\zeta_{\eta,0,j}\bboxtimes \zeta_{\eta,1,j}$,
  where
   $G_{\eta}: {\mathfrak M}_{\Gamma_0^o}\times {\mathfrak M}_{\Gamma_1^o}
              \rightarrow {\mathfrak M}_{g,n}$ is the natural morphism
   between the related moduli stack of nodal curves.
 Then
 $$
  \Psi_{(g,n;\beta)}^{W^{(3)}_t}(\alpha(t),\zeta)\;
   =\; \sum_{\eta\in\overline{\Omega}_{(g,n;\beta)}^{\cal L}}\;
        \frac{{\mathbf m}(\eta)}{|\Eq(\eta)|}\,
          \sum_{j\in K_{\eta}}\,
           \left[ \Psi_{\Gamma_0}^{Y_0^{\rel}}(j_0^{\ast}\alpha(0),
                    \zeta_{\eta,0,j})\,
                  \bullet\,
                  \Psi_{\Gamma_1}^{Y_1^{\rel}}(j_1^{\ast}\alpha(0),
                    \zeta_{\eta,1,j})
           \right]_0\,,
 $$
 where
  $j_i:Y_i\rightarrow W^{(3)}_0$ is the inclusion map,
  $$
   \Psi_{\Gamma_i}^{Y_i^{rel}}(j_i^{\alpha}(0),\zeta_{\eta,i,j})
   = {\mathbf q}_{i\,\ast}\left(
       \ev^{\ast}(j_i^{\ast}\alpha(0))
      \cdot \pi_{\Gamma_i}^{\ast}(\zeta_{\eta,i,j})
      \cdot  [{\mathfrak M}({\mathfrak Y}_i^{\rel},\Gamma_i)]^{\virt}
                          \right) \in H_{\ast}(D^r)\;, \;\; i=0,1\,.
  $$

 In cycle form,
 $$
  [{\mathfrak M}({\mathfrak W}^{(3)}_0,(g,n;\beta))]^{\virt}\;
   =\; \sum_{\eta\in\overline{\Omega}_{(g,n;\beta)}^{\cal L}}\;
         \frac{{\mathbf m}(\eta)}{|\Eq(\eta)|}\,
             \Phi_{\eta\ast}\Delta^!
               \left(
                [{\mathfrak M}({\mathfrak Y}_0^{\rel},\Gamma_0)]^{\virt}
                \times
                [{\mathfrak M}({\mathfrak Y}_1^{\rel},\Gamma_1)]^{\virt}
               \right)\,,
 $$
 where
 ${\mathfrak M}({\mathfrak W}^{(3)}_0,(g,n;\beta))$
  is the fiber of ${\mathfrak M}({\mathfrak W}^{(3)},(g,n;\beta))$
  over ${\mathbf 0}\in{\Bbb A}^1$,
 $\Delta^!$ is the Gysin map associated to the diagonal map
 $\Delta:D^r\rightarrow D^r\times D^r$ for the relevant $D^r$
  in each summand.
} 

\bigskip

We want to make things as intrinsic to $Y$ as possible so that
 we can appropriately combine the Gromov-Witten invariants of
 $Y_0\cup_D Y_1$ defined from [Li2] to a quantity that is
 justifiable to be called a Gromov-Witten invariant of $Y$.
In particular, we want to remove the possible ${\cal L}$-dependence
 (so far on {\it both} sides of the equation) in the above gluing
 formula.
Let us now deal with this issue.

\bigskip

\begin{flushleft}
{\bf The ${\cal L}$-(in)dependence of $\Omega^{\cal L}_{(g,k;\beta)}$
     and ${\mathfrak M}({\mathfrak W}^{(3)}_0,(g,n;\beta))$.}
\end{flushleft}
Recall from Sec.~1 the intermediate families over ${\Bbb A}^1$
 that occur in the semi-stable reduction of $W/{\Bbb A}^1\,$:
 $$
  \begin{array}{cccccccl}
   W^{(3)}  & \stackrel{\varphi_3}{\longrightarrow}  & W^{(2)}
            & \stackrel{\varphi_2}{\longrightarrow}  & W^{(1)}
            & \stackrel{\varphi_1}{\longrightarrow}  & W^{(0)} := W \\
   \mbox{\scriptsize $\pi_3$}\downarrow\hspace{3ex}
    && \mbox{\scriptsize $\pi_2$}\downarrow\hspace{3ex}
    && \hspace{2ex}\downarrow\mbox{\scriptsize $\pi_1$}
    && \hspace{6ex}\downarrow\mbox{\scriptsize $\pi_0=\pi$}\\[.6ex]
   {\Bbb A}^1   & =  & {\Bbb A}^1   & \stackrel{\alpha}{\longrightarrow}
    & {\Bbb A}^1     & =  & {\Bbb A}^1   &\,.
  \end{array}
 $$
Denote
 the exceptional divisor of $\varphi_1$ on $W^{(1)}$ by $E_1$
  ($\simeq{\Bbb P}^3$) and
recall that the exceptional divisor of $\varphi_3$ on $W^{(3)}$
 is $Y_1$.
Recall [C-H] or [L-Y: Remark 3.1].
Let ${\cal L}_0$ be a sufficiently very ample line bundle on
 $W^{(0)}=W$.
Then ${\cal L}_1:=(\varphi_1^{\ast}{\cal L}_0)(-E_1)$ is very ample
 on $W^{(1)}$.
The pull-back ${\cal L}_2:=\varphi_2^{\ast}{\cal L}_1$ on $W^{(2)}$
 is $\pi_2$-ample.
Thus there is an open subset $U$ of ${\Bbb A}^1$, containing
 ${\mathbf 0}$, such that ${\cal L}_2$ is ample on $\pi_2^{-1}(U)$.
By removing finitely many fibers of the families/${\Bbb A}^1$
 in the above diagram and with an abuse of notation that we denote
  $\pi_2^{-1}(U)$ over $U$ still by $W^{(2)}/{\Bbb A}^1$
  (and since it is only a neighborhood of $Y$ in $W$ that matters),
 we will say that ${\cal L}_2$ is ample on $W^{(2)}$.
By taking $k>\!>0$, ${\cal L}_2^{\otimes k}$ becomes sufficiently
 very ample on $W^{(2)}$ and
its pull-back with a twist
 ${\cal L}:= (\varphi_3^{\ast}{\cal L}_2^{\otimes k})(-Y_1)$
 becomes very ample on $W^{(3)}$.

\bigskip

\noindent
{\bf Lemma 2.2
  [$\,{\cal L}$-independence of $\Omega^{\,{\cal L}}_{(g,n;\beta)}\,$].}
{\it
 Let
  ${\cal L}$ be a very ample line bundle on $W^{(3)}$ as constructed
   above and
  $d$ be the ${\cal L}$-degree of $\beta\in H_2(Y;{\Bbb Z})$.
 Then $\Omega^{\,{\cal L}}_{(g,n;\beta)}$ that appears in the
  degeneration formula in Lemma 2.1
  is independent of ${\cal L}$.
} 

\bigskip

\noindent
{\it Proof.}
The proof follows the same reasoning as in the proof of
 [L-Y: Lemma 3.2].
Recall the blow-up resolution $p_0:Y_0\rightarrow Y$ with exceptional
 divisor identical to $D$ ($=Y_0\cap Y_1$)
 and $p_1$ that sends $Y_1$ to the conifold singularity of $Y$.
It follows from a careful chasing through the sequence of
 pull-backs and twists in the construction of ${\cal L}$ that
 $$
  {\cal L}^{(0)}:={\cal L}|_{Y_0}\;
    \simeq\; (p_0^{\ast}{\cal L}_0|_Y)^{\otimes k}(-(k+1)D)
   \hspace{1em}\mbox{and}\hspace{1em}
  {\cal L}^{(1)}:={\cal L}|_{Y_1}\;
              \simeq\; {\cal O}_{Y_1}(-(k+1)D)\,.
 $$
Recall the generator $\gamma_2$ of $H_2(Y_1;{\Bbb Z})$.
Let $\gamma_{2,1}$ and $\gamma_{2,2}$ be the curve classes on $Y_0$
 from the two rulings of $D$. Then they generate the semi-group
 of curve classes in the relative Mori cone $\NE(p_0)$.
Observe that $D\cdot \gamma_{2,i}=-1$ on $Y_0$ for $i=1,2$
 while $D\cdot\gamma_2=+1$ on $Y_1$.
Let $\NE(Y)_{\beta}$ be the set of all the curve classes
 in $\NE(Y)_{\scriptsizeBbb Z}$ that represent $\beta$
 in $H_2(Y;{\Bbb Z})$.
Then, for each $\gamma\in\NE(Y)_{\beta}$,
 $p_{0\ast}^{-1}(\gamma)$ in $\NE(Y_0)_{\scriptsizeBbb Z}$ is of the form
  $\widetilde{\gamma}
   +{\Bbb Z}_{\ge 0}\gamma_{2,1}+{\Bbb Z}_{\ge 0}\gamma_{2,2}$
  for a unique $\widetilde{\gamma}\in\NE(Y_0)_{\scriptsizeBbb Z}$
  determined by $\gamma$.

Now let $\eta=(\Gamma_0,\Gamma_1,I)\in \Omega_{(g,n;\beta)}^{\cal L}$.
Then the pairs
 $(b(\Gamma_0),b(\Gamma_1))$ are characterized by the conditions:
 $$
  p_{0\ast} b(\Gamma_0)\; =\;\beta
  \hspace{1em}\mbox{and}\hspace{1em}
  {\cal L}^{(0)}\cdot b(\Gamma_0)+{\cal L}^{(1)}\cdot b(\Gamma_1)\;
  =\; {\cal L}\cdot \beta\; =\; d\,.
 $$
Solving it explicitly as in [L-Y: proof of Lemma 3.2],
 one concludes that

{\footnotesize
 $$
  \Omega^{\,{\cal L}}_{(g,n;\beta)}\;
  =\;
     \coprod_{\gamma\in NE(Y)_{\beta}}\;
         \left\{
          \begin{array}{l}
            \eta=(\Gamma_0,\Gamma_1,I) \\
             \mbox{admissible} \\
             \mbox{triple}     \\
             \mbox{for $Y_0\cup_D Y_1$}
          \end{array}
          \left|
            \begin{array}{l}
              \bullet\hspace{1ex}
                  b(\Gamma_0)=\widetilde{\gamma}
                              +l_{0,1}\gamma_{2,1}+l_{0,2}\gamma_{2,2},\,
                  b(\Gamma_1)=l_1\gamma_2\,, \\
              \hspace{2ex} (l_{0,1}+l_{0,2})+l_1
                             =D\cdot \widetilde{\gamma}\,,\;
                           l_{0,1},\,l_{0,2},\,l_1
                                    \in{\footnotesizeBbb Z}_{\ge 0}\,;
                                                                 \\[.6ex]
              \bullet\hspace{1ex}
               g(\eta) = g\,,\;
               k_1+k_2 = n\,; \\[.6ex]
              \bullet\hspace{1ex}
               \sum_{i}\mu_{0,i}\,=\,l_1\,;\\[.6ex]
              \bullet\hspace{1ex}
               I\subset \{1,\,\ldots,\, n\}\,,\; |I|=k_1\,.
            \end{array}
          \right.
         \right\}\;
 $$
} 

\noindent
This set is indeed independent of ${\cal L}$.

\noindent\hspace{14cm}$\Box$

\bigskip

\noindent
We shall denote the ${\cal L}$-independent set of admissible triples
 worked out explicitly in the end of the proof above by
 $\Omega_{(g,n;\beta)}$.

\bigskip

The proof of the above lemma implies also that, with this choice of
 ${\cal L}$ on $W^{(3)}/{\Bbb A}^1$, the potentially
 ${\cal L}$-dependent stack
 ${\mathfrak M}({\mathfrak W}^{(3)},(g,n;\beta))$
 can be re-defined without referring to ${\cal L}$ at all and, hence,
 depends only on $W$ and $(g,k;\beta)$.
Its fiber ${\mathfrak M}({\mathfrak W}^{(3)}_0,(g,n;\beta))$
 over ${\mathbf 0}\in{\Bbb A}^1$ is thus also ${\cal L}$-independent.
Adding this ${\cal L}$-independence into the statements in
 Lemma 2.1, 
 we obtain:

\bigskip

\noindent
{\bf Theorem 2.3 [intrinsic to $W/{\Bbb A}^1$ and $\beta$].}
(Cf.\ [L-Y: Theorem 3.3].)
{\it
 The degeneration/gluing formulas, one in the numerical form and
  the other in the equivalent cycle form,
 in Lemma 2.1 are independent of ${\cal L}$ and are intrinsic to
  $W$ and $(g,n;\beta)$ when ${\cal L}$ is chosen as above:
 $($the cycle form omitted$\,)$
 $$
  \Psi_{(g,n;\beta)}^{W^{(3)}_t}(\alpha(t),\zeta)\;
   =\; \sum_{\eta\in\overline{\Omega}_{(g,n;\beta)}}\;
        \frac{{\mathbf m}(\eta)}{|\Eq(\eta)|}\,
          \sum_{j\in K_{\eta}}\,
           \left[ \Psi_{\Gamma_0}^{Y_0^{\rel}}(j_0^{\ast}\alpha(0),
                    \zeta_{\eta,0,j})\,
                  \bullet\,
                  \Psi_{\Gamma_1}^{Y_1^{\rel}}(j_1^{\ast}\alpha(0),
                    \zeta_{\eta,1,j})
           \right]_0\,.
 $$
} 

\bigskip

\noindent
{\it Remark 2.4
 $[\,W$-dependence of
     ${\mathfrak M}({\mathfrak W}^{(3)}_0,(g,n;\beta))\,]$.}
From the details in [Li2: Sec.~4.4], the dependence of the stack
 ${\mathfrak M}({\mathfrak W}^{(3)}_0,(g,n;\beta))$
 on the neighborhood of $Y$ in $W$, if any, can be only mild.
The fact that the morphism
 $$
  \coprod_{\eta\in\Omega_{(g,n;\beta)}}\,\Phi_{\eta}\;:\;
   \coprod_{\eta\in\Omega_{(g,n;\beta)}}\,
   {\mathfrak M}({\mathfrak Y}_0^{\rel},\Gamma_0)\times_{D^r}\,
   {\mathfrak M}({\mathfrak Y}_1^{\rel},\Gamma_1)\;
   \longrightarrow\; {\mathfrak M}({\mathfrak W}^{(3)}_0,(g,n;\beta))
 $$
 is natural and surjective ([Li1: Sec.~4] and [Li2: Sec.~4.4]) on
 the underlying topological space of points with the Zariski topology
implies that different choices of $W$ for a given $Y$
 give stacks ${\mathfrak M}({\mathfrak W}^{(3)}_0,(g,n;\beta))$
 of homeomorphic underlying topologies, (cf.\ [L-MB: Chapter~5]).
Since each $\Phi_{\eta}$ is finite \'{e}tale of degree
 a $W$-independent combinatorial factor $|\Eq(\eta)|$ to its image
 ${\mathfrak M}(
     {\mathfrak Y}_0^{\rel}\sqcup{\mathfrak Y}_1^{\rel},\eta)$
 ([Li1: Sec.~4]),
as long as these stacks have the same relative multiplicity to
 the corresponding image of $\Phi_{\eta}$ -
 which is proved to be true and this relative multiplicity of
  the component of ${\mathfrak M}({\mathfrak W}^{(3)}_0,(g,n;\beta))$
  labelled by $\eta$ is given by ${\mathbf m}(\eta)$ irrelevant to $W$
  ([Li2: Sec.~4.4]) -
any detailed difference will not affect the resulting value of
 Gromov-Witten invariants of $Y$ to be extracted from these stacks
 for a given $(g,n;\beta)$.
(Such indifference is implicit in [Li2] in order to define 
  Gromov-Witten invariants of the degenerate fiber of a family
  that is family-independent and
 indeed the degeneration formula implies this indifference as well.)

\bigskip

\bigskip

\section{Extracting GW invariants of a conifold from pairs.}

The results in Sec.~1 ans Sec.~2 together with the
 constant-under-deformation requirement for any good definition
 of Gromov-Witten invariants of conifolds
implies the following route to extract Gromov-Witten invariants
 of the conifold $Y$ from the pairs $(Y_0,D)$ and $(Y_1,D)$,
(in this section, by ``Gromov-Witten invariants of $Y$", we mean
 the invariants of $Y$ that would be defined via the intersection
  theory on the moduli stack $\overline{\cal M}_{g,n}(Y,\beta\,)$
  of stable maps to $Y$, 
 should $\overline{\cal M}_{g,n}(Y,\beta\,)$ exists and admits
  a perfect obstruction theory):

{\footnotesize
$$
 (g,n;\beta\,)\; \Rightarrow\;
 \left\{\,
  \left(
   {\mathfrak M}({\mathfrak Y}_0^{\rel},\Gamma_0)\,,\,
   {\mathfrak M}({\mathfrak Y}_1^{\rel},\Gamma_1)
  \right)\,
 \right\}_{\eta=(\Gamma_0,\,\Gamma_1,\,I)
                      \in\overline{\Omega}_{(g,n;\beta\,)}}\;
 \stackrel{\mbox{\parbox{6.8em}{\tiny
   gluing formula \newline from Theorem 2.3\newline}}}{\Rightarrow}\;
 \mbox{\parbox{8em}{\scriptsize
   Gromov-Witten \newline invariants of $Y$ \newline
   as would be defined  \newline
   via $\overline{\cal M}_{g,n}(Y,\beta\,)\,$.}}
$$
} 

\bigskip

There is one last ingredient, though, that we have not yet discussed:
 when the number $n$ of marked points on the prestable curve is
 non-zero, the definition of Gromov-Witten invariants of $Y$
 involves choices of cycles on $Y$.
In this case, we need to know whether we can choose canonically
 (perhaps up to some equivalence relation) a set of cycles on
 $Y_0\cup_D Y_1$ to feed into the expression in the definition
 of relative Gromov-Witten invariants for $(Y_0,D)$ and $(Y_1,D)$
 in the gluing formula.
Let us now turn to this last issue.

Suppose that $Y=W_0$ is in a degeneration $W/{\Bbb A}^1$.
Let $\tau_i:A_{\ast}(Y_i)\rightarrow A_{\ast-1}(D)$, $i=0,\,1$,
 be the group homomorphism defined by taking intersection with $D$.
Let $Z$ be a cycle in $W$ flat over ${\Bbb A}^1$,
 then $(Z\cdot Y_0)\cdot D = - (Z\cdot Y_1)\cdot D$.
Since $Y_0\cdot D = - Y_1\cdot D$ on $W$, the equality 
 $\tau_0(Z\cdot Y_0)=\tau_1(Z\cdot Y_1)$ must hold.
This motivates the following definition:

\bigskip 

\noindent
{\bf Definition 3.1 [pre-deformable class on $Y_0\cup_D Y_1$].}
{\rm
 Let $\xi_i\in A_{\ast}(Y_i)$ , $i=0,\,1$. We say that
  $(\xi_0,\xi_1)$ gives a {\it pre-deformable class}
  (i.e.\ $\xi_0+\xi_1$) on $Y_0\cup_D Y_1$
 if $\tau_0(\xi_0)=\tau_1(\xi_1)$ in $A_{\ast}(D)$.
 Such $(\xi_0,\xi_1)$ will be called an {\it admissible pair}
  of cycle classes with respect to the gluing $Y_0\cup_D Y_1$.
} 

\bigskip

\noindent
Readers may note that the image of a pre-deformable map from
 a (not necessarily connected) prestable curve to $Y_0\cup_D Y_1$
 as defined in [Li1: Sec.~2; in particular, Definition 2.9]
 gives an example of a pre-deformable class on $Y_0\cup_D Y_1$.

Recall ([Fu]) that the rational map $\varsigma_0:Y\dasharrow Y_0$
 (the inverse of $p_0:Y_0\rightarrow Y$) defines a group
 homomorphism
 $\varsigma_{0\ast}:A_{\ast}(Y)\rightarrow A_{\ast}(Y_0)$
 by taking intersection product with the closure of the graph of
 $\varsigma_0$ in $Y\times Y_0$.

\bigskip

\noindent
{\bf Definition 3.2 [admissible class on $Y$].} {\rm
 A class $\xi\in A_{\ast}(Y)$ will be called {\it admissible}
  if $\varsigma_{0\ast}(\xi)$ in $A_{\ast}(Y_0)$ admits
  a unique extension to an admissible pair
  $(\varsigma_{0\ast}(\xi),\xi^{\prime})$ with respect to
  $Y_0\cup_D Y_1$.
 Denote the set of all such classes in $A_{\ast}(Y)$ by
  $A_{\ast}^{ad}(Y)$.
} 

\bigskip

Note the $\tau_1:A_{\ast}(Y_1)\rightarrow A_{\ast}(D)$
 is injective on $A_{\ge 1}(Y_1)$ with $\tau_1(A_{\ge 1}(Y_1))$
 spanned by $[\pt]$, [$(1,1)$-curve], and $[D]$.
Consequently, effective classes of constant dimension $\xi_1$
 in $A_{\ge 1}(Y_0)$ whose intersection with $D$ lies in
 this subspace has a unique extension to an admissible pair
 $(\xi_1,\xi_2)$ with respect to $Y_0\cup_D Y_1$.
The image of such $\xi_1$ on $Y$ will span $A_{\ast}^{ad}(Y)$.
The constant-under-deformation requirement and multi-linearity
 in the axioms of Gromov-Witten invariants for any variety imply
 the following:

\bigskip

\noindent
{\bf Corollary 3.3 [extraction of GW invariants of $Y$ from pairs].}
{\it
 The Gromov-Witten invariants of $Y$ that involve only classes
  in $A_{\ast}^{ad}(Y)$ can be extracted from the relative
  Gromov-Witten invariants of pair $(Y_0,D)$ and $(Y_1,D)$
  by the gluing formula.
} 

\bigskip

In other words, 
let
 $\alpha\in (A_{\ast}^{ad}(Y)_{\scriptsizeBbb Q})^{\times n}$ and
 $\zeta\in A_{\ast}({\mathfrak M}_{g,n})$.
Denote by $\Psi_{(g,n;\beta)}^Y(\alpha,\zeta)$
  the Gromov-Witten invariant of $Y$ associated to the
  topological type $(g,n;\beta\,)$ and classes $\alpha$ and $\zeta$
 that is defined via {\it any}$\,$ standard procedure,
 i.e.\ a construction of $\overline{\cal M}_{g,n}(Y,\beta)$
  with a perfect obstruction theory, ..., etc.\ in such a way
  that the Gromov-Witten axioms and the constant-under-deformation
   property are satisfied and
  that when the construction is applied to smooth variety,
   it recovers the Gromov-Witten invariants from the equivalent
    construction of [B-M] and [L-T].
Let
 $(\alpha_0,\alpha_1)\in
  A_{\ast}(Y_0)^{\times n}\times A_{\ast}(Y_1)^{\times n}$
 be the tuple from the admissible pairs with respect to $Y_0\cup_D Y_1$
 of the class in the entries of $\alpha$.
For
 \mbox{$\eta=(\Gamma_0,\,\Gamma_1,\,I)\in\overline{\Omega}_{(g,n;\beta)}$},
assume that
 $G_{\eta}^{\ast}(\zeta)
   =\sum_{j\in K_{\eta}}\zeta_{\eta,0,j}\bboxtimes \zeta_{\eta,1,j}$,
 where
  \mbox{$G_{\eta}:
         {\mathfrak M}_{\Gamma_0^o}\times {\mathfrak M}_{\Gamma_1^o}
             \rightarrow {\mathfrak M}_{g,n}$}
  is the natural morphism between the related moduli stack of nodal
  curves, cf.\ [Li2: Sec.~0].
Then
$$
 \Psi_{(g,n;\beta)}^Y(\alpha,\zeta)\;
  =\; \sum_{\eta\in\overline{\Omega}_{(g,n;\beta)}}\;
       \frac{{\mathbf m}(\eta)}{|\Eq(\eta)|}\,
         \sum_{j\in K_{\eta}}\,
          \left[ \Psi_{\Gamma_0}^{Y_0^{\rel}}(\alpha_0,
                   \zeta_{\eta,0,j})\,
                 \bullet\,
                 \Psi_{\Gamma_1}^{Y_1^{\rel}}(\alpha_1,
                   \zeta_{\eta,1,j})
          \right]_0\,.
$$
While the left-hand side of the identity remains to be truly
 constructed
 (though one may turn this identity to a {\it definition}
  of Gromov-Witten invariants of $Y$ for $\alpha$ admissible
  if one wishes),
 the right-hand side is completely determined by
 and canonically/intrinsically associated to the data:
 $Y$, $(g,n;\beta\,)$, $\alpha$, and $\zeta$.

Note that when $Y$ is a Calabi-Yau conifold (i.e.\ it arises from
 a degeneration of Calabi-Yau $3$-folds), the expected dimension
 of the would-be moduli stack $\overline{\cal M}_{g,0}(Y,\beta\,)$
 is zero.
In this case, no issue of choice of cycles is involved. 

\bigskip

\noindent
{\bf Corollary 3.4 [Calabi-Yau conifold].} {\it
 When $Y$ is a Calabi-Yau conifold, all the Gromov-Witten invariants
  of $\,Y$ associated to the topological type $(g,0\,;\,\beta\,)$
  can be extracted from the relative Gromov-Witten invariants of
  pairs $(Y_0,D)$ and $(Y_1,D)$ by the gluing formula.
} 

\bigskip

\noindent
{\it Remark 3.5 $[\,$general conifold$\,]$.} {\rm
 For a general conifold $Y$ with $\Sing(Y)=\{\,p_1,\,\cdots\,\}$,
  the discussions, expressions, and formulas in Sec.~2 and Sec.~3
  generalize immediately
  by appropriately do the following replacement:
  $D$ in $Y_0$ by $\coprod_i D_i$ in $Y_0$,
  $(Y_1,D)$ by isomorphic copies of smooth quadric pairs
   $(Y_i,D_i)\simeq (Y_1,D)$,
  summand involving the index $1$ of the pair $(Y_1,D)$ by
   the summation $\sum_i\,$, ..., etc., in the expressions involved,
   where $i=1,\,\ldots.$ corresponds to each conifold singularity
   $p_i$ of $Y$.
} 

\bigskip

\noindent
{\it Remark 3.6 $[\,$a general picture$\,]$.}
For a general irreducible singular variety $Y$, the result in
 the current work suggests the following picture whose details remain to be studied:
 an assignment to $Y$ a collection of pairs:
  a principal pair $(Y_0,D_0=\cup_iD_{0,i})$,
   where $Y_0\rightarrow Y$ is a resolution of $Y$ whose exceptional
    locus $D_0$ is a divisor with only simple normal crossing (snc)
    singularities, and
  a set of auxiliary smooth (smooth variety)-(snc divisor) pairs
   $(Y_j,D_j)$ that are related to the germs of the singularities
   of $Y$ only,
  and a gluing diagram $\Gamma$ that encodes the gluing all the these
   pairs 
  so that
 (cf.\ {\sc Figure}~3-1) 
 \begin{itemize}
  \item[(1)]
  the relative Gromov-Witten invariants of $(Y_0,D_0)$
   is obtained from a combination of those for the smooth pairs
   $(Y_0,D_{0,i})$, $(D_{0,i}, D_{0,i}\cap \cup_{i^{\prime}\ne i}D_{0,j\i^{\prime}})$,
    ..., etc,
  and similarly for $(Y_j,D_j)$
  (cf.\ a ``quantum inclusion-exclusion principle");

  \item[(2)]
   the Gromov-Witten invariants of $Y$ is then obtained from
    a combination of the relative Gromov-Witten invariants of
    $(Y_0,D_0)$ and $(Y_j,D_j)$'s;

  \item[(3)]
  the combination rules in Item (1) and Item (2) above are
   functorial/universal in the same sense that
   Jun Li's degeneration formula is functorial/universal
   (i.e.\ independent of details of $Y$ away from the singularities,
          curve classes chosen, cycle classes chosen, ... etc.).
 \end{itemize}

\begin{figure}[htbp]
 \setcaption{{\sc Figure} 3-1.
  \baselineskip 14pt
  A gluing construction of Gromov-Witten invariants of a singular variety $Y$.
  Since Gromov-Witten invariants arise from A-model topological
   string theory on the physics side [Wi3], such gluing construction
   should ring closely with the Atiyah-Segal's formulation of
   topological quntum field theory ([At] and [Se];
   see also [E-G-H] and [B-P]).
 } 
\centerline{\psfig{figure=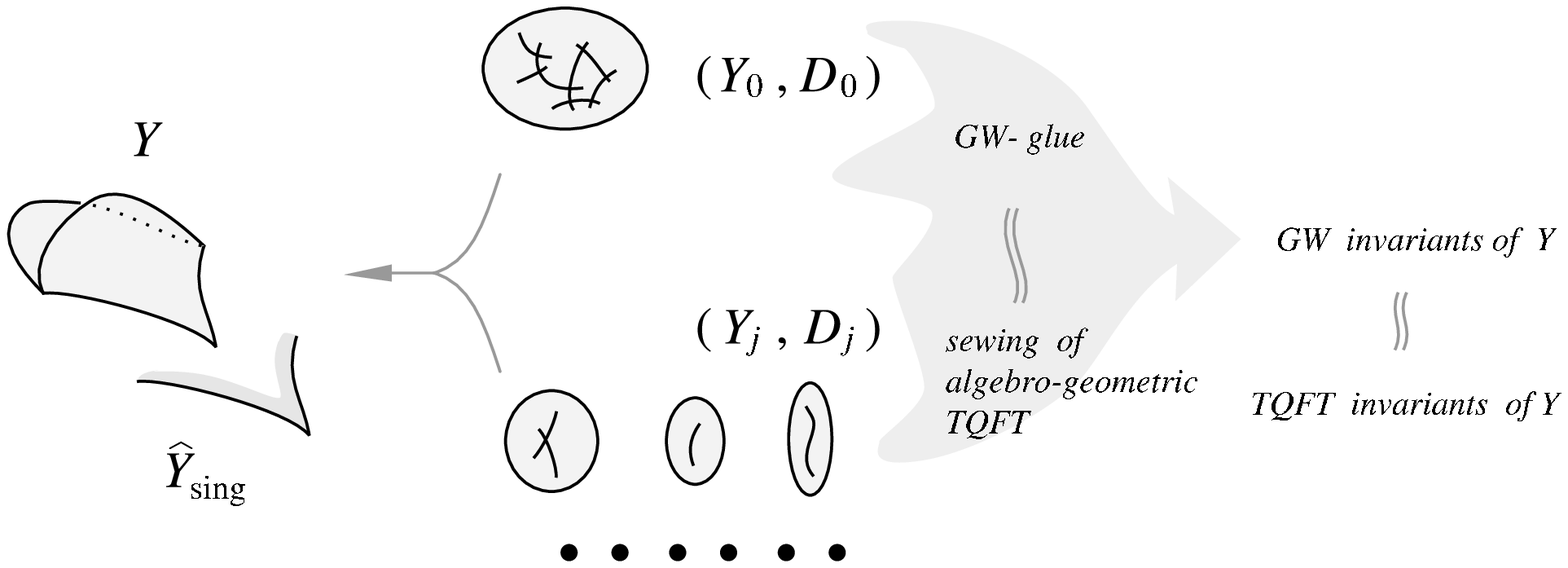,width=13cm,caption=}}
\end{figure}

\vspace{6em}

\end{document}